
\input amstex 

\def\tri{\therosteritem}
 
\def\bcc#1{\Bbb C^{#1}}
\def\pp#1{\Bbb P^{#1}}
\def\ppp{\Bbb P }\def\bccc{\Bbb C }

\def\ci{\Cal I}
\def\co{\Cal O}

\def\hd{, \hdots ,}

\def\ot{\!\otimes\!}

\def\pii{\partial II}

\def\ra{\rightarrow}
\def\rtwo{\frac r 2}

\def\tdim{\text{dim}}\def\tcodim{\text{codim}}
\def\thom{\text{Hom}}
\def\tker{\text{ker}}
\def\tmax{\text{ max }}

\def\trank{\text{rank}}

\def\ww{\wedge}

\def\xsm{X_{sm}}

\def\lx#1{x_{#1}}
\def\lh#1{H_{#1}}

\define\intprod{\mathbin{\hbox{\vrule height .5pt width 3.5pt depth 0pt %
        \vrule height 6pt width .5pt depth 0pt}}}
\documentstyle{amsppt}
\magnification = 1200
\hsize =15truecm
\hcorrection{.5truein}
\baselineskip =18truept
\vsize =22truecm
\NoBlackBoxes

\NoRunningHeads

\topmatter
\title
On symmetric degeneracy loci, spaces of symmetric matrices of constant
rank and dual varieties
\endtitle

\leftheadtext{Bo Ilic and J.M. Landsberg}
\rightheadtext{On symmetric degeneracy loci \ldots }
\author
 Bo Ilic${}^*$ and J.M. Landsberg${}^{**}$
\endauthor

\date { October 9, 1996}\enddate

\address {Department of Mathematics, UCLA,
Los Angeles, CA 90095}\endaddress
\email {ilic\@math.ucla.edu} \endemail
\address{Department of Mathematics,
Columbia University, New York,  NY 10027}
\endaddress
\curraddr Math\'ematique,
Laboratoire Emile Picard,
Universit\'e Paul Sabatier,
31062 Toulouse CEDEX,
France  \endcurraddr
\email {jml\@picard.ups-tlse.fr }
\endemail

\thanks {* Supported by an NSERC postdoctoral fellowship.}
\endthanks
\thanks {** Supported by NSF grant DMS-9303704.}
\endthanks
\endkeywords
\subjclass{14, 15, 53}\endsubjclass
\abstract{ 
  Let $X$ be a nonsingular simply connected projective
variety of dimension $m$, $E$ a rank $n$ vector bundle on $X$, and $L$
a line bundle on $X$. Suppose that $S^2(E^{*}) \otimes L$ is an ample 
vector bundle and that there is a constant even rank $r \ge 2$
symmetric bundle map $E \to E^{*} \otimes L$. We prove 
that $m \le n-r$.
We use this result to
solve the constant rank problem for symmetric matrices, proving
that the maximal dimension of a linear subspace of the space
of $m\times m$ symmetric  matrices such that each nonzero element has even rank
$r \ge 2$ is $m-r+1$.  We explain how this result relates to 
the study of dual varieties in projective geometry
and give some applications and examples. }
\endabstract

\endtopmatter

\document

\heading  \S 0. Introduction \endheading

Let $V$ be a complex vector space of dimension $m$. Consider a linear 
subspace $A \subset S^2 V^*$. We can think of $A$ as
a vector space of symmetric $m \times m$ matrices or
equivalently, as a symmetric $m \times m$ matrix whose entries are linear 
forms on $A$. We say that $A$ has {\it constant rank} $r$ if every nonzero
element of $A$ has rank $r$.  We ask two basic questions about constant
rank subspaces $A$: how large can such subspaces be and what
invariants can we attach to them? These questions come up in a number of
contexts: linear algebra, vector bundles on projective space,
symmetric degeneracy loci and dual varieties.

In linear algebra, an important motivation for this work is the
classical theory, due to Kronecker and Weierstrass, giving a 
normal form for singular pencils of matrices or quadratic forms.
That is, if $A \subset \thom( \bcc{m}, \bcc{n})$ or $A \subset 
S^2V^*$ is a linear subspace such that each element of $A$ has rank
$\le r$ (i.e. $A$ is of {\it bounded rank} $r$) and $\tdim(A) = 2$,
then a convenient normal form can be given 
from which geometric properties of the pencil are easy to read off.
An excellent exposition of both cases is given by F. R. Gantmacher 
\cite{G, Chp. XII} or, for the case of quadratic forms one can
consult Hodge and Pedoe \cite{HP, Chp. VIII, Sec. 10}.
If $A \subset S^2 V^*$
is  a pencil of constant rank $r$, the normal form shows
that in fact $r$ must be even. R. Meshulam proved a generalization of
this result without using the normal form for the pencil \cite{M}. We
also prove it in two different ways in Theorem 2.15.

The classification problem when $\tdim(A) \ge 3$ is more difficult although
several interesting results are known.
R. Westwick, building on work of J. Sylvester, showed that if a linear 
subspace $A \subset \thom( \bcc{m}, \bcc{n})$ has constant rank $r$, then
$\tdim(A) \le m + n - 2r + 1$ \cite{W}, \cite{S}. 
H. Flanders proved  that if $r \le m \le n$, and
$A \subset \thom( \bcc{m}, \bcc{n})$ is a linear space of bounded
rank $r$ then $\tdim(A) \le rn$ \cite{F} (see also \cite{M2}). There is
also a classification of such subspaces  of near maximal
dimension with the most recent results due to L.B. Beasley \cite{B}.
R. Meshulam found  bounds for symmetric and skew symmetric linear families
analogous to Flanders' bound for general linear families mentioned
above \cite{M3}.
Another approach
is to try to classify subspaces $A \subset \thom( \bcc{m}, \bcc{n})$
of bounded rank $r$ for small values of $r$. This was accomplished
first by M.D. Atkinson for $r \le 3$ and then reproved  
by D. Eisenbud and J. Harris \cite{A}, \cite{EH}.

Our initial interest in linear spaces $A \subset S^2 V^*$ of constant 
rank $r$ arose because such spaces can be constructed from
smooth varieties having degenerate dual varieties. Recall that if
$X^n \subset \pp{N}$ is a nonsingular projective variety of dimension $n$
then the {\it dual variety} $X^* \subset \pp{N*}$ is the union of all 
hyperplanes $H$ such such that $H \cap X$ is singular, i.e. such that $H$
is tangent to $X$. A dimension count shows that
one expects $X^*$ to be a hypersurface and it is interesting to study when
this fails to occur. To that end,  define the {\it defect} of $X$, 
$\delta$ by $\delta := N-1 - \tdim(X^*)$ and  assume that $X^*$ is
degenerate, i.e. $\delta \ge 1$. Let $H$ be a smooth point of $X^*$
and let $|II_{X^*, H}|  \subset \ppp (S^2T^*_H X^*)$ 
denote the linear system of quadrics generated by the second fundamental
form of $X^*$ at $H$ (see e.g. \cite{Lan}). We prove:

\proclaim{Theorem 3.4} 
$|II_{X^*,H}|$ is a linear system of quadrics of projective dimension
$\delta$ and constant rank $n-\delta$.
\endproclaim

Thus, given a smooth variety with degenerate dual variety, 
Theorem 3.4 constructs  a 
$\delta +1$ dimensional linear subspace $A \subset S^2 \bcc{N-1-\delta}$
of constant rank $n-\delta$. This constant rank space of quadrics can
also be viewed as arising from a result of L. Ein \cite{E, Theorem 2.2}.
We also prove an explicit inversion formula that allows one to recover
$II_{X^*,H}$ from the second fundamental form and cubic form of
any smooth point on $X$ that $H$ is tangent to (Theorem 3.9).
Also, Griffiths and Harris proved that linear spaces 
$A \subset S^2 V^*$ of bounded rank $r$ can be constructed from
not necessarily smooth varieties which have degenerate duals \cite{GH2, 3.5}.
We give details and complete references in \S 3.

Our main linear algebraic result is:

\proclaim{Theorem 2.16} If $r$ is even and $\ge 2$ then
$$\tmax \{ \tdim(A) \mid A \subset S^2V^* \text{ is of constant rank } r \}
 = m-r+1.$$
\endproclaim

If $X_r \subset \ppp(S^2 V^*)$ is the projective variety of 
rank $\le r$ matrices, then an elementary dimension count shows that
$\tcodim(X_{r-1}, X_r) = m - r + 1$
which explains why one expects this bound. As
noted above, the result when $r$ is odd is classical; see Theorem 2.15. 
As an application of Theorem 2.16 and Theorem 3.4
we give a new proof of F. Zak's result that for a smooth variety $X$
with dual variety $X^*$, $\tdim(X^*) \ge \tdim(X)$. Moreover, Theorem 2.15
and Theorem 3.4 give a proof of the Landman parity theorem:
if $X$ is a nonsingular projective variety with degenerate dual variety then
$n-\delta$ is even (notation as above).

We prove Theorem 2.16 by first proving a more general step-wise result
on symmetric degeneracy loci:

\proclaim{Theorem 1.2} Let $X$ be a nonsingular simply connected projective
variety of dimension $m$, $E$ a rank $n$ vector bundle on $X$, and $L$
a line bundle on $X$. Suppose that $S^2(E^{*}) \otimes L$ is an ample 
vector bundle and that there is a constant even rank $r \ge 2$
symmetric bundle map $E \to E^{*} \otimes L$. Then $m \le n-r$.
\endproclaim

To $A \subset S^2V^*$ we associate a short exact sequence of vector bundles on
$\ppp(A)$:
$$0 \to K \to \co_{\ppp(A)}^{m} \to E \to 0$$
where $E$ has rank $r$ and
$E \cong E^*(1)$ and then investigate restrictions on $E$ and its 
chern classes in 2.4-2.6. We give examples of symmetric
(and skew symmetric) linear spaces $A$ and compute the associated exact
sequences in 2.9-2.14.

We also survey some of the results of the linear algebraists,
most recently due to Westwick, sometimes providing alternate proofs 
and references and tying them in with the algebraic geometry literature.

\subheading{Acknowledgements}
It is a pleasure to thank Henry Pinkham and Mark Green for helpful
discussions and encouragement. We are especially indebted
to Robert Lazarsfeld for his suggestions towards Theorem 1.2.

\heading \S 1. Even rank symmetric degeneracy loci \endheading

For the basic properties of ample vector bundles and the Lefschetz
theorem used in this section we refer to \cite{Laz}. If $W$ is
a variety, we use the notation $b^i(W) := \tdim (H^i(W,\bccc))$.
The following is a well known proposition that will be used in the proof of
the next theorem: 
\proclaim{Proposition 1.1} Let $Q^{r-2} \subset \pp{r-1}$ be a nonsingular
quadric hypersurface. Then.
\roster
\item $b^i(Q) = 1$ if $0 \le i \le 2(r-2)$, $i$ is even and $i \ne r-2$.
\item $b^{r-2}(Q) = 2$ if $r$ is even.
\item $b^i(Q) = 0$ for all other cases.
\endroster
\endproclaim
 
\demo{Proof}
Except for $b^{r-2}(Q)$, all the cohomology of $Q$ can be computed 
by using the Lefschetz theorem, Poincar\'{e}
duality and the cohomology of projective space. 
But $e(Q)$, the topological Euler characteristic of $Q$
is just the degree of the top chern class $\int_Q c_{r-2}(T_Q)$
which is $r$ if $r$ is even and $r-1$ if $r$ is
odd. (The chern class 
 computation is done for nonsingular surfaces in $\pp{3}$ in
\cite{GH, pg 601}. The computation for nonsingular hypersurfaces is
 analogous).
This determines $b^{r-2}(Q)$.
\qed
\enddemo

\proclaim{Theorem 1.2} Let $X$ be a nonsingular simply connected projective
variety of dimension $m$, $E$ a rank $n$ vector bundle on $X$, and $L$
a line bundle on $X$. Suppose that $S^2(E^{*}) \otimes L$ is an ample 
vector bundle and that there is a constant even rank $r \ge 2$
symmetric bundle map $E \to E^{*} \otimes L$. Then $m \le n-r$.
\endproclaim

\demo{Proof}
Consider the projective bundle  $\pi: \ppp(E) \to X$. We first show that 
there is a subvariety $Y \subset \ppp(E)$  given as the zero locus
of a section $t \in H^0(\ppp(E), \co_{\ppp(E)}(2) \otimes \pi^{*}L)$
such that the fiber of $\pi \vert_Y : Y \to X$ over each $x \in X$
is a rank $r$ quadric hypersurface in $\ppp(E(x))$.

Indeed, there is a natural map 
$\pi^* \pi_{*} \co_{\ppp(E)}(2) \to \co_{\ppp(E)}(2)$. Since
$\pi_{*} \co_{\ppp(E)}(2) \cong S^2(E^*)$ this gives after tensoring with
$\pi^{*}L$ a map $u: \pi^* ( S^2(E^*) \otimes L) \to
 \co_{\ppp(E)}(2) \otimes \pi^* L$.
On the other hand, the symmetric bundle map $E \to E^{*} \otimes L$ defines
a section $s \in H^0(X,S^2(E^*) \otimes L)$ which pulls back to give
a section $\pi^*s \in H^0(\ppp(E), \pi^* ( S^2(E^*) \otimes L))$. Let 
$t = u \circ \pi^*s \in  H^0(\ppp(E), \co_{\ppp(E)}(2) \otimes \pi^{*}L)$.

The section $t$ at $[v] \in \ppp(E(x))$ is the linear map
$\lambda \to s(x)(v,v) \lambda$. and so vanishes iff
$s(x)(v,v) = 0$. But by hypothesis, this defines a rank r quadric hypersurface
in $\ppp(E(x))$.

\proclaim{Claim 1.3} $\co_{\ppp(E)}(2) \otimes \pi^{*}L$ is ample.
\endproclaim

\demo{Proof}
First, recall that if $E$ is any vector bundle and $L$ any line bundle
then  $\ppp(E \otimes L) \cong \ppp(E)$ and 
$\co_{\ppp(E \otimes L^*)}(1) \cong \co_{\ppp(E)}(1) \otimes \pi^*L$.

Then since $S^2(E^*) \otimes L$ is ample on $X$,
$\co_{\ppp(S^2(E) \otimes L^*)}(1) \cong \co_{\ppp(S^2 E)}(1)
 \otimes \sigma^* L$ is ample on $\ppp(S^2 E)$ where
$\sigma: \ppp(S^2(E)) \to X$ is the natural projection.
Now the second Veronese gives an inclusion $i: \ppp(E) \subset 
\ppp(S^2 E)$ such that $\pi = i \circ \sigma$. Then
$i^* \co_{\ppp(S^2 E)}(1) \cong \co_{\ppp(E)}(2)$. So
$i^* ( \co_{\ppp(S^2 E)}(1) \otimes \sigma^* L) \cong \co_{\ppp(E)}(2)
\otimes \pi^* L$ is ample on $\ppp(E)$ as required.
\qed
\enddemo

Thus, $\ppp(E) \setminus Y$ is an affine variety and has
the homotopy type of a CW complex of (real) dimension $\le m+n -1$.
Thus $H_i(\ppp(E) \setminus Y) = 0$ for $i >= m+n$.

Let $K$ denote the kernel and $F$ denote the image of the map
$E \to E^* \otimes L$. Then since the map has constant rank $r$,
$K$ and $F$ are vector bundles on $X$ of ranks $n-r$ and $r$ respectively.
Since the map
is symmetric, there is a symmetric isomorphism $F \cong F^* \otimes L$.
We have a natural map
$p: \ppp(E) \setminus \ppp(K) \to \ppp(F)$ given
on the fibers by linear projection  centered at $\ppp(K(x)) \subset 
\ppp(E(x))$. $p$ is a $\bcc{n-r}$- fiber bundle.

The isomorphism $F \cong F^* \otimes L$ determines a
 hypersurface $Z \subset \ppp(F)$
such that if $\rho: \ppp(F) \to X$ is the natural projection then
the fiber of  $\rho \vert_Z : Z \to X$ over each $x \in X$
is a smooth quadric in $\ppp(F(x))$.

Now $\ppp(K) \subset Y$ and so  $p$ restricts to a
$\bcc{n-r}$- fiber bundle map $\ppp(E) \setminus Y \to \ppp(F) \setminus Z$.
Thus, $H_i(\ppp(E) \setminus Y) \cong H_i(\ppp(F) \setminus Z)$ and by
Lefschetz duality this is isomorphic to $H^{2(m+r-1) -i}(\ppp(F), Z)$.
Hence $H^i(\ppp(F), Z) = 0$ for $i \le 2r + m - n -2$ and then using the
long exact sequence of the pair $(\ppp(F), Z)$ we conclude that:
$H^i(Z) \cong H^i(\ppp(F))$ for $i \le 2r + m - n -3$.

Now $b^i(\ppp(F)) = b^i(X \times \pp{r-1})$.
So if $f^{p,q} = b^p(X) b^q(\pp{r-1})$ then by the Kunneth formula,
$b^i(\ppp(F)) = \sum_{p+q = i} f^{p,q}$.

By Deligne's theorem \cite{GH, pg 466}, the Leray spectral sequence for 
$\rho \vert_Z : Z \to X$
degenerates at the $E_2$ term. Since $X$ is simply connected there is no
monodromy so 
$(R^i \rho \vert_Z) (\bccc) = \text{ the constant sheaf } H^i(Q,\bccc)$
 where
$Q \subset \pp{r-1}$ is a smooth quadric. Thus,
$E_2^{p,q} = H^p(X) \otimes H^q(Q)$.
 Let $e_2^{p,q} = \text{dim} (E_2^{p,q})$.
Then $b^i(Z) = \sum_{p+q=i} e_2^{p,q}$.

Let $g^{p,q} = e_2^{p,q} - f^{p,q}$. Thus by the above,
 $\sum_{p+q=i} g^{p,q} = 0$ for $i \le 2r+m+n-3$.
 On the other hand, since $b^{r-2}(Q) = 2$ (this is where we use the fact 
that $r$ is even), $\sum_{p+q = r-2} g^{p,q} = b^0(X) = 1$.
Thus $2r+m-n-3 < r-2$ and so $m \le n-r$ as required.
\qed
\enddemo

The previous theorem can be viewed as a version for symmetric maps
of the following theorem of R. Lazarsfeld:
\proclaim{Theorem 1.4} \cite{Laz, Theorem 2.2, pg. 41}
 Let $X$ be a projective variety of dimension $m$. Let
$E$ and $F$ be vector bundles on $X$ of ranks $e$ and $f$ respectively.
Suppose that $E^* \otimes F$ is ample and that there is a constant rank
$r$ vector bundle map $E \to F$. Then $m \le e + f - 2r$.
\endproclaim

Lazarsfeld used this theorem to give a step-wise proof of the
non-emptyness of degeneracy loci which was originally proved 
by Fulton and Lazarsfeld:
\proclaim{Theorem 1.5} \cite{Laz, Theorem 2.1, pg. 40}
 Let $X$ be a projective variety of dimension $m$. Let
$E$ and $F$ be vector bundles on $X$ of ranks $e$ and $f$ respectively.
Suppose that $E^* \otimes F$ is ample and let $\phi:E \to F$ be
a vector bundle map. If $m \ge (e-r)(f-r)$ then
$X_r(\phi) = \{ x \in X \mid \trank(\phi(x)) \le r \}$ is non-empty.
\endproclaim

There are analogues of the above theorem for symmetric maps and
skew symmetric maps due to Harris and Tu \cite{T}, \cite{HT}.
However, in the case that
$r$ is odd, the proven result is not as sharp as the result conjectured.

\subheading{Question 1.6} Is Theorem 1.2 true without the assumptions that
$X$ is nonsingular and simply connected and that the rank $r$ is even?  
An affirmative answer would give a step-wise proof of the non-emptiness
of symmetric degeneracy loci. 

\heading \S 2. Linear spaces of matrices satisfying rank conditions \endheading

Let $V=\bcc{m}$, $W= \bcc{n}$ and let $A \subset \thom(V,W)$ 
(or $A\subset S^2V^*$ or $A\subset\Lambda^2V^*)$   be a linear subspace.
$A$ is said to be of {\it bounded rank $r$} if  for
all $q\in A$, $\trank (q)\leq r$,
of {\it   rank bounded below by $r$} if  for
all nonzero $q\in A$, $\trank (q)\geq r$,
 and of
{\it constant rank $r$} if  for
all nonzero $q\in A$, $\trank (q)= r$.

\subheading{Definition 2.1}
$$
\align
 &l(r,m,n)=\tmax\{ \tdim (A) \ | \ A\subset V\ot W
\text{ is of constant rank }r\}\\
&\underline l(r,m,n)=\tmax\{ \tdim (A) \ | \ A\subset V\ot W
\text{ is of bounded rank  }r\}\\
&\overline l(r,m,n)=\tmax\{ \tdim (A) \ | \ A\subset V\ot W
\text{ is of   rank bounded below by }r\}\\
\endalign
$$

Similarly, define $c(r,m)$, $\underline c(r,m)$, and $\overline c(r,m)$ in
the symmetric case and $\lambda(r,m)$, $\underline \lambda(r,m)$, and
$\overline \lambda(r,m)$ in the skew symmetric case.

\subheading{2.2 Vector bundles and maps 
associated to linear spaces of matrices}

For the basic concepts and definitions of vector bundles on projective
space used in this part (e.g. chern classes, the splitting principle
and uniform bundles) we refer to \cite{OSS}. We say that a vector bundle
$E$ is {\it free} if it splits as a direct sum of line bundles.
Given a linear subspace $A \subset \thom (V,W)$ with $\tdim(A) = l+1$ we 
can associate a vector bundle map 
$\psi: V \otimes \co_{\ppp(A)}(-1) \to W \otimes \co_{\ppp(A)}$ on
$\ppp(A)$ as follows:
at $[q] \in \ppp(A)$ the fiber of $\co_{\ppp(A)}(-1)$ is 
$\lambda q, \lambda \in \bccc$ and 
$\psi(v \otimes \lambda q) =\lambda \cdot q(v)$. 
Tensoring by $\co_{\ppp(A)}(1)$
we get a map $\phi: \co_{\pp{l}}^m \to \co_{\pp{l}}^n(1)$. Alternatively,
choosing a basis for $V$ and $W$ and a basis $x_0, \dots x_l$ for $A^*$,
$A$ can be viewed as an $n \times m$ matrix of linear forms in the $x_i$ and
such a matrix gives $\phi$.

If $A$ has constant rank $r$ then the kernel $K$, cokernel $N$ and image $E$
are vector bundles of rank $m-r$, $n-r$ and $r$ respectively and determine
short exact sequences:
$ 0 \to K \to \co_{\pp{l}}^m \to E \to 0$
and $0 \to E \to \co_{\pp{l}}^n(1) \to N \to 0$.

Let $h = c_1( \co_{\pp{l}})$. If $F$ is any vector bundle on $\pp{l}$,
we can write $c_i(F) = f_i h^i$ for some $f_i \in {\Bbb Z}$. Thus
$c(F) = 1 + c_1(F) + \dots + c_l(F) = 1 + f_1 h + \dots f_l h^l$.

The following theorem is due to R. Westwick \cite{W}. The above construction
and the use of chern classes to obtain a weaker results along the
lines of the next theorem are due to J. Sylvester \cite{S}. 
\proclaim{Theorem 2.3} \cite{W} Suppose $2 \le r \le m \le n$. Then
\roster
\item $l(r,m,n) \le m + n - 2r +1$
\item $l(r,m,n) = n - r + 1$ if $n-r+1$ does not divide $(m-1)! / (r-1)!$
\item $l(r,r+1,2r-1) = r +1$
\endroster
\endproclaim

\demo{Proof}
We'll sketch Westwick's proof of \tri1 and the $\le$ direction of \tri2.
(The $\ge$ direction is classical; see Proposition 2.10 below).
From the above exact sequences, $c(K)c(E) = 1$ and $c(E)c(N) = (1+h)^n$.
Thus, $c(K)(1+h)^n = c(N)$. If $n-r+1 \le i \le l$ then $c_i(N) = 0$ and
looking at the coefficient of $h^i$ we get
$\sum_{j = 0}^{m-r} \binom{n}{i-j} k_j = 0$ where 
$k_j h^j = c_j(K)$ and 
we use the convention that 
$\binom{n}{j} = 0$ if $j<0$ or $j > n$.  The coefficient matrix of this
collection of linear equations is 
$M = ( \binom{n}{i-j})_{0 \le j \le m-r, \  n-r+1 \le i \le l}$.
If $l = m + n -2r +1$ then this is a square invertible matrix with
determinant $\prod_{j=0}^{m-r} j!$. Thus $k_0 = 0$ which is a contradiction
since $k_0 = 1$. This proves \tri1.
Westwick refers to a privately published manuscript of Muir and
 Metzler as a reference for evaluating
this determinant however one can also refer to e.g.  \cite{ACGH, pg. 93-95}.
\tri2 follows directly from considering the linear equation with $i = n-r+1$.

We remark that \tri1 is also a direct consequence of the vector
bundle construction above and Lazarsfeld's result, Theorem 1.4. (This was 
also previously observed by R. Meshulam). 
\qed
\enddemo

If $A \subset S^2V^*$ (resp. $A \subset \Lambda^2 V^*$) we similarly get
a symmetric (resp. skew symmetric) map
$\phi: \co_{\pp{l}}^m \cong  V \otimes \co_{\ppp(A)} \to
  V^* \otimes \co_{\ppp(A)}(1)
 \cong \co_{\pp{l}}^m (1)$ 
which can be viewed as given by an $m \times m$ symmetric (resp. skew
symmetric) matrix of linear forms in the $x_i$.

If $A$ has constant rank $r$, the symmetry or skew symmetry
 implies that $N \cong K^*(1)$ and
$E \cong E^*(1)$. This is because dualizing and then twisting the
map $\phi: \co_{\pp{l}}^m \to \co_{\pp{l}}^m (1)$ by $\co_{\pp{l}}(1)$
we get the same map $\phi$ in the symmetric case and $-\phi$ in the skew
symmetric case. The two short exact sequences above thus reduce to the
single sequence $0 \to K \to \co_{\pp{l}}^m \to E \to 0$
with $E \cong E^*(1)$.

Conversely, given a surjection $\co_{\pp{l}}^m \to E \to 0$ where $E$ is
a vector bundle of rank $r$ satisfying $E \cong E^*(1)$ we can dualize and
twist the surjection by  $\co_{\pp{l}}(1)$ to obtain:
$\co_{\pp{l}}^m \to E \cong  E^*(1) \to \co_{\pp{l}}^m(1)$. The above 
composition then is given by an $m \times m$ matrix of linear forms  and leads
to a dimension $l+1$ linear space $A \subset \thom(\bcc{m}, \bcc{m})$ of
constant rank $r$. Note that the subspace $A$ need not be symmetric or
skew symmetric. For example, it could be the direct sum of symmetric and
skew symmetric subspaces.

\proclaim{Proposition 2.4}
 Let $A \subset S^2 V^*$ or $A \subset \Lambda^2V^*$ be a constant
rank $r$ subspace of dimension $\ge 2$. 
 Then $E$ is a uniform vector bundle of splitting type
$\co_{\pp 1}^{\frac r2}\oplus \Cal O^{\frac r2}_{\pp 1}(1)$. In particular,
$r$ is even.
\endproclaim

\demo{Proof}
Since $E \cong E^*(1)$, $E$ must be uniform of splitting type
$(a_1\hd a_{\frac r2},b_1\hd b_{\rtwo})$ where $b_i=1-a_i$ (here one
may need to reorder). But since $E$ is globally generated,
all the factors in its splitting type must be nonnegative.
\qed\enddemo

Using the isomorphism $E \cong E^*(1)$ we get linear relations on the
chern classes of $E$. Specifically, suppose that $r \ge l$, which is the
only interesting case by Proposition 2.6. If $c_i(E) = e_i h^i$ 
then for $0 \le i \le l$, 
$$e_i = \sum_{j=0}^{i} \binom{r-j}{i-j}(-1)^j e_j.$$
For $i$ odd, this 
expresses $e_i$ as a linear combination of $e_j$'s with $j < i$. Thus
there are $\ge [l/2]$ linearly independent relations and it is not
too hard to see that there are in fact exactly $[l/2]$. 
If $i = 1$ we get $e_1 = r/2$ (which gives another proof that
$r$ is even) and if $i=3$, using $e_1=r/2$ we get
$(r-2)e_2 - 2e_3 = r(r-1)(r-2)/12$.

\proclaim{Proposition 2.5}
 For $1 \le i \le n$, $0 \le e_i \le e_1^i = (r/2) ^i$.
Thus $\{ (c_1(E), \dots, c_l(E)) \mid E$  is a globally generated 
vector bundle of rank $r$ on $\pp{l}$ and $E \cong E^*(1) \}$
is a finite set. 
\endproclaim

\demo{Proof}
This follows by the discussion in \cite{DPS, pg 317}.
\qed
\enddemo

\proclaim{Proposition 2.6} If $r \le l$ then 
$E \cong \co_{\pp{l}}^{r/2} \oplus \co_{\pp{l}}^{r/2}(1)$
 or $r=l=2$ and $E \cong T_{\pp{2}}(-1)$.
\endproclaim

\demo{Proof}
If $E$ is free then $E \cong \co_{\pp{l}}^{r/2} \oplus \co_{\pp{l}}^{r/2}(1)$
by Propostition 2.4. 
If $r < l$ then since $E$ is uniform by \cite{OSS, Theorem
3.2.3, pg 55} $E$ is free.
If $r=l$ and $E$ is not free then  by the discussion
 on \cite{OSS, pg 71}, $E \cong T_{\pp{l}}(a)$ or $\Omega_{\pp{l}}(a)$.
Since $c_1(E) = l / 2$, taking first chern classes we conclude that
$E  \cong T_{\pp{2}}(-1)$ or $\Omega_{\pp{2}}(2)$; but these
last two bundles are isomorphic.
\qed
\enddemo

\subheading{2.7 Linear spaces with rank bounded below by $r$}

Let $X_r = \{ [q] \in \ppp(\thom(V,W)) \mid \trank (q) \le r \}$ 
and similarly for $X_r \subset \ppp(S^2 V^*)$ and 
$X_r \subset \ppp(\Lambda^2 V^*)$. Then an elementary dimension count
(see e.g \cite{ACGH, pg. 67 \& pg. 101}) shows that:

\roster
\item $\tcodim (X_r \subset \ppp(\thom(V,W))) = (m-r)(n-r)$
\item $\tcodim (X_r \subset \ppp(S^2 V^*)) = \binom{m-r+1}{2}$
\item $\tcodim (X_r \subset \ppp(\Lambda^2 V^*))= \binom{m-r}{2}$ ($r$ even)
\endroster
   
By Bezout's theorem, this implies the essentially classical:
\proclaim{Proposition 2.8}
\roster
\item $\overline l(r,m,n) = (m-r)(n-r)$.
\item $\overline c(r,m) = \binom{m-r+1}2$.
\item $\overline \lambda(r,m) = \binom{m-r}{2}$ ($r$ even).
\endroster
\endproclaim


\subheading{2.9 Some examples of symmetric and skew-symmetric constant
 rank linear systems}

Given a linear subspace $A \subset \thom(V,W)$, we can consider $A$ as
a linear subspace of $S^2(V^* \oplus W)$ via the natural inclusion
$\thom(V,W) = V^* \otimes W \subset S^2(V^* \oplus W)$. To distinguish 
these two different linear embeddings of $A$,
we'll denote $A \subset  S^2(V^* \oplus W)$ by $B$.
If $A$ has constant rank (resp. bounded rank,
 resp. rank bounded below by) $r$ then $B$ is a linear space of quadrics
of constant rank (resp. bounded rank,
 resp. rank bounded below by) $2r$ . We say that $B$ is a {\it doubling}
 of $A$.
In matrices, systems formed by doubling look like:
$$
\left\{ \pmatrix 0&a\\ {}^ta& 0\endpmatrix  \mid  a\in A\right\}.
$$
Similarly, we can double $A \subset \thom(V,W)$ to get 
$B \subset \Lambda^2(V^* \oplus W)$.

The following proposition is essentially classical:
\proclaim{Proposition 2.10} 
\roster
\item If $0 <r \le m \le n$ then $l(r,m,n) \ge n - r + 1$.
\item If $r$ is even and $\ge 2$
then $c(r,m) \ge m-r+1$ and $\lambda(r,m) \ge m-r+1$.
\endroster
\endproclaim

\demo{Proof}
$X_{r-1} \subset \ppp(\thom(\bcc{r}, \bcc{n}))$ has codimension $n-r+1$. Thus
we can find an $n-r+1$ dimensional
 subspace $A \subset \thom(\bcc{r}, \bcc{n})$ of
constant rank $r$ from which \tri1 follows.
For \tri2, note that doubling $A$ produces an $n-r+1$ dimensional subspace
$B \subset S^2 \bcc{r+n}$
(or $\Lambda^2 \bcc{r+n}$) of constant rank $2r$. 
\qed 
\enddemo

\subheading{2.11 Example of linear systems with $E$ free}
For the  $m-r+1$ dimensional subspace $B$
of $S^2 \bcc{m}$ (or $\Lambda^2 \bcc{m}$) of constant even rank $r$ produced
in the above proof the associated vector bundle 
$E \cong \co_{\pp{m-r}}^{r/2} \oplus \co_{\pp{m-r}}^{r/2}(1)$.
This follows for the symmetric case since in matrices
$$
B = \pmatrix 0 & A \\ {}^tA & 0 \endpmatrix =
\pmatrix 0 & I_{r/2} \\ {}^t A & 0 \endpmatrix
\pmatrix I_{r/2} & 0 \\ 0 & A \endpmatrix
$$
and so the image of the induced map $\phi: \co_{\pp{m-r}}^m \to \co_{\pp{m-r}}^m(1)$
is clearly $\co_{\pp{m-r}}^{r/2} \oplus \co_{\pp{m-r}}^{r/2}(1)$.

\subheading{2.12 Examples of linear systems with $E$ not free
arising from the Euler sequence}

Given any vector space $V$ of dimension $m$, for $1 \le k \le m-1$ there
is a natural inclusion
$$
\align &V\ra \thom (\Lambda^k V, \Lambda ^{k+1}V)  \\
& v\to (\phi_v: \alpha \to v \ww \alpha)\endalign
$$
Now if $v \ne 0$, $\tker(\phi_v) = v \ww \Lambda^{k-1}V$ and so this gives an 
$m$-dimensional linear space of linear maps of constant rank
$\binom {m}k -\binom {m-1} {k-1} = \binom{m-1}{k}$.
 This observation is due to Atkinson and Westwick \cite{AW, pg 233}.

Now suppose that $\tdim(V) = 2a+1$ and fix a volume form on $V$ i.e. an 
isomorphism $\Lambda^{2a+1}V \cong \bccc$ so that we can identify
$\Lambda^{a}V^* \cong \Lambda^{a+1}V$. Then $V \subset \thom(\Lambda^{a}V,
\Lambda^{a}V^*)$ and  we can think of $\phi_v$ as a bilinear form on
$\Lambda^{a}V$. Then 
 $\phi_v(\alpha, \beta) = v \ww \alpha \ww \beta = (-1)^a v \ww \beta \ww
\alpha = (-1)^a \phi_v(\beta, \alpha)$. Thus, if $a$ is even then in fact
$V \subset S^2(\Lambda^{a}V^*)$ and if  $a$ is odd then
$V \subset \Lambda^2(\Lambda^{a}V^*)$.
So, we have constructed a $2a+1$ dimensional
linear space of $\binom{2a+1}{a} \times \binom{2a+1}{a}$
matrices (symmetric for $a$ even and skew symmetric for $a$ odd)
of rank exactly $\binom{2a}{a}$.

From the vector bundle point of view, at $[v] \in \ppp(V)$, the map
$\Lambda^a V \to \Lambda^a V^*$ factors:
$\Lambda^a V \to \Lambda^a (V/v) \to \Lambda^a V^*$. Recall that the fiber
of $T_{\ppp(V)}(-1)$ at $[v]$ is $V/v$ so we conclude that $E$ is just
$\Lambda^{a} (T_{\ppp(V)}(-1))$. Thus the associated short exact sequence is
just the $a$-th exterior power of the Euler sequence i.e.
$$
0 \to \co_{\pp{2a}}(-1) \otimes \Lambda^{a-1} (T_{\pp{2a}}(-1)) \to
\Lambda^a \co_{\pp{2a}}^{2a+1} \to \Lambda^{a} (T_{\pp{2a}}(-1)) \to 0
$$
In particular, $E$ is not free.

The $a=1$ case yields the $3\times 3$ skew
symmetric matrices; doubled it is the linear
system of quadrics corresponding to the second fundamental
form of the Grassmanian $G(2,5)$ (which is self-dual) (see \S 3),
which is also the space of quadrics vanishing on
the Segre, $\pp 1\times \pp 2 \subset \pp{5}$.
The $a=2$ case (which is already symmetric) is the
linear system of quadrics corresponding to the second fundamental
form of the ten dimensional spinor variety (which is 
also self-dual), which is also the space of quadrics vanishing on
the Grassmanian, $G(2,5)$. (One can prove these statements via a direct
coordinate computation).
  
\subheading{2.13 Westwick's example} 
What follows is an intrinsic construction of a 3 dimensional linear space
of $2a+1 \times 2a+1$ skew symmetric matrices of constant rank $2a$ where
$a \ge 1$. For $a=1$ this is the $3 \times 3$ skew symmetric matrices.
This example  is given in coordinates in Sylvester's paper 
\cite{S, pg 4} where it is attributed to Westwick.
Let $V$ be a $2a+1$ dimensional vector space. 
Then $\ppp := \ppp(\Lambda^2 V^*) = X_{2a}$
 and $X_{2a-1}= X_{2a-2}$ since a skew
symmetric matrix can not have odd rank. Now $\tcodim( X_{2a-2}, \ppp) = 3$ so
any $\ppp^2 \subset \ppp \setminus X_{2a-2}$ will work.
An easy chern class computation gives that the associated short exact 
sequence is:
$$
0 \to \co_{\pp{2}}(-a) \to \co_{\pp{2}}^{2a+1} \to E \to 0.
$$

\subheading{Example 2.14}
Consider the special case when $l=3$ and $m - r = 2$. Such
an example can't arise from a linear space of symmetric matrices of constant
rank by Theorem 2.16. Then  $0 = c_3(K) = -e_3 + 2e_1 e_2 -e_1^3$.
Using the  linear relations in the  $e_i$ from the discussion following
Proposition 2.4  
it follows that $c(K) = 1 + (-r/2)h + (r(r+1)/12)h^2$ and 
$c(E) = 1 + (r/2)h + (r(2r-1)/12)h^2 + (r^2(r-2)/24)h^3$ and so $r=8$ gives
the smallest rank possibility: $c(E) = 1 + 4h + 10h^2 + 16h^3$. In fact, 
such an example, given by a skew symmetric family $A$ was given by
Westwick \cite{W2, pg 168}.

\proclaim{Theorem 2.15} If $r$ is odd then $c(r,m) = 1$.
\endproclaim

This theorem follows immediately from Proposition 2.4
or alternatively, from the first chern class computation in the
discussion after Proposition 2.4. However, as mentioned in the
introduction, it was known
classically as a consequence of the Kronecker-Weierstrass theory giving
a normal form for pencils of symmetric matrices of bounded rank $r$
\cite{G, Chp. XII}, \cite{HP, Chp. VIII, Sec. 10}.
There is also a generalization by Roy Meshulam \cite{M}.

\proclaim{Theorem 2.16} If $r$ is even and $\ge 2$ then $c(r,m) = m - r + 1$.
\endproclaim

\demo{Proof}
By Proposition 2.10 we only need to prove the $\le$ inequality.
Given $A \subset S^2V^*$ of constant rank r 
with $\tdim(A) = l + 1$ and $\tdim(V) = m$
as above we obtain a symmetric constant even rank $r$ map 
$\co_{\pp{l}}^m \to \co_{\pp{l}}^m(1)$. Now $\pp{l}$ is simply connected
and $(S^2 \co_{\pp{l}}^{m}) \otimes \co_{\pp{l}}(1)$ is ample and so
by Theorem 1.2, $l <= m-r$ as required.
\enddemo

\heading \S 3. Applications to dual varieties \endheading

For the basic background on dual varieties 
for this section see \cite{E} and the references
there. Let $X^n \subset \pp{N}$ be a nonsingular projective variety and let
$X^* \subset \pp{N*}$ be the dual variety of $X$. Define the defect
of $X$, $\delta$, by $\delta := N-1 - \tdim(X^*)$. Assume that $\delta \ge 1$
and
let $H$ be a  smooth point of $X^*$. Then $H$ can
be considered as a hyperplane in $\pp{N}$ and the contact locus
$L = L_H = \{ x \in X \mid \tilde T_x X \subset H \} $ is a $\delta$ 
dimensional (projective) linear space.
 Let $N_{L/X}$ be the normal bundle to $L$ in $X$.
In his work on dual varieties Lawrence Ein proved that there is an isomorphism 
$N_{L/X} \cong N_{L/X}^*(1)$ \cite{E, Theorem 2.2}.
In fact, it follows from Ein's proof that this isomorphism is symmetric.

It is convenient to reformulate Ein's result as
follows in order to  construct vector spaces of symmetric matrices of
constant rank:
\proclaim{Lemma 3.1} A nonsingular projective variety
$X^n \subset \pp{N}$  with defect $\delta > 0$ determines a
$\delta +1$ dimensional linear subspace $A \subset S^2 \bcc{N-1-\delta}$
of constant rank $n-\delta$ with associated short exact sequence:
$$
0 \to N_{X/H}^*(1) \vert_L \to \co_{\pp{\delta}}^{N-1-\delta}
\to N_{L/X} \to 0.
$$  
\endproclaim

\demo{Proof}
Start with the exact sequence
$$
0 \to N_{L/X} \to N_{L/H} \to N_{X/H} \vert_L \to 0
$$
which holds since $H$ is tangent to $X$ along $L$. 
Then dualize, twist by $\co_{\pp{\delta}}(1)$ and use the fact that
$N_{L/H} \cong \co_{\pp{\delta}}^{N-1-\delta}(1)$ and $N_{L/X} \cong N_{L/X}^*(1)$
to obtain the short exact sequence in the statement of the lemma. Thus 
since $N_{L/X} \cong N_{L/X}^*(1)$ is a symmetric isomorphism, 
we can obtain the constant rank linear subspace of symmetric matrices as 
described in 2.2.
\qed
\enddemo
 
The above lemma, along with Theorem 2.16 gives a new proof of a theorem
of F. Zak:
\proclaim{Theorem 3.2} \cite{Z, Cor. 10, pg. 39}
Let $X^n \subset \pp{N}$ be a nonsingular projective variety.
Then $\tdim(X^*) \ge \tdim(X)$. 
\endproclaim

Let $x\in X$ be a smooth point and let $II_{X,x}\in S^2T^*_xX\ot N_xX$
denote the projective second fundamental form of $X$ at $x$
(see e.g. \cite{Lan}).
Let $|II_{X,x}|=\ppp (II^*(N^*_xX)) \subset \ppp (S^2T^*_xX)$ 
denote the linear system of quadrics it generates. 

The following theorem is a result proved (but not stated explicitly) by
Griffiths and Harris \cite{GH2, 3.5}. It shows that linear spaces of
symmetric matrices of bounded rank arise from varieties 
$X$ with degenerate duals. Note that the result does not assume that
$X$ is nonsingular. 
\proclaim{Theorem 3.3} \cite{Lan, 5.2}
Let $X^n\subset\pp{N}$ be a  variety with defect $\delta \ge 1$ and let
$x\in X$ be a general point. Then
$|II_{X,x}|$ is a linear system of quadrics of projective dimension
$\delta$ and bounded rank $n - \delta$.
\endproclaim

Another interpretation of Lemma 3.1 is:  
\proclaim{Theorem 3.4} 
Let $X^n \subset \pp{N}$ be a smooth variety with defect $\delta \ge 1$ 
and let $H\in X^*$ be any smooth point. Then
$|II_{X^*,H}|$ is a linear system of quadrics of projective dimension
$\delta$ and constant rank $n-\delta$.
\endproclaim

We remark that Theorem 3.4 provides a nice geometric way of seeing the
linear system of quadrics of constant rank whose study is the basis
of this paper. In fact, this observation was the beginning of our 
investigation of this problem and only later did we see the connection
with the prior work of Ein. 

We present our original proof of (3.4) as 
we derive a slightly more refined formula which shows how
a singular quadric in $|II_{X^*,H}|$ exactly corresponds
to a singular point of $X$.

Write $\pp N=\ppp V$ and 
let $\ci\subset X\times X^*$ denote
the standard  incidence correspondence.
Denote the smooth points of a variety $Y$ by
$Y_{sm}$. Let $H\in \xsm^*$ and
let $x\in L_H\cap \xsm$.
Let $<,>: V\times V^*\ra \Bbb C$ denote the canonical pairing.
We let $\hat T_xX\subset V$ denote the
cone over the embedded tangent space and in general a hat on
an object in projective space denotes the corresponding cone. 
Note that $<,>$ descends to   a    
  nondegenerate pairing between
$r$ dimensional vector spaces:
$$
(\hat T_xX/\hat T_xL_H)\times (\hat T_HX^*/\hat  T_HL_x)\ra \Bbb C \tag 3.5
$$
which we continue to denote by $<,>$.
(To see why, recall that $\hat T_xL_H= N^*_HX^*(1)$.)
Let $(\lx t,\lh t)$ be the lifting of a curve in
$\ci$ to $V\times V^*$. Note that
$<\lx 0',\lh 0>=0$, and
$ <\lx 0,\lh 0'>=0$, where prime denotes derivative.
Differentiating $<\lx 0,\lh 0'>$, 
  we have
$$
<\lx 0',\lh 0'>=- 
<\lx 0,\lh 0''>.\tag 3.6
$$
Now $\lx 0'\in \hat T_xX, \lh 0'\in \hat T_HX^*$ and
up to twists, 
$<\lx 0,\lh 0''>= n_x\intprod II_{X^*,H}(W,W)$
where $W\in T_HX^*$ corresponds to $\lh 0'$
and $n_x\in N^*_HX^*$ corresponds to $\lx 0$.
We polarize (3.6) to conclude that the quadratic form
$n_x\intprod II_{X^*,H}$ is of rank $r$. If we assume
that
  $L_H\cap X_{sing}=\emptyset$, then
$x$ corresponds to an arbitrary point of $\ppp N^*_HX^*$
and we conclude:

\proclaim{Theorem 3.7} Let $X^n\subset\pp N$ be a variety
with dual variety 
$X^*$ having defect $\delta = n-r$.
Say there exists an $H\in X^*_{sm}$   such that
$H$ is only tangent to smooth points of $X$. Then
$|II_{X^*,H}|$ is a $\delta$ dimensional system of quadrics
of constant rank $r$ on $T_HX^* \cong \bcc{N-\delta-1}$.
\endproclaim

\smallpagebreak

Note that dual to (3.6) we have
$
<\lx 0',\lh 0'>=- 
<\lx 0'',\lh 0 >, 
$
which together with (3.6)
implies the following result:

\proclaim{Theorem 3.8} Let $x\in \xsm$ and let $H\in X^*_{sm}$ be
tangent to $x$. Then
$$
q_H=Q_x
$$
where $q_H\subset \ppp (T_xX/T_xL_H)$ is the smooth quadric hypersurface
corresponding to evaluating $II_{X,x}$ at a vector representing
$H$ and projecting to the quotient,  
$Q_x\subset\ppp ( T_HX^*/T_HL_x)$
is the corresponding dual object, and $T_xX/T_xL_H$ is canonically
identified
with $(T_HX^*/T_HL_x)^*$ up to twists by (3.5) (and both spaces are
  identified
with their duals by virtue of the quadrics on them).
\endproclaim

To see the connection with 
(\cite{Ein}, 2.2), note that
$
T_xX/T_xL_H= N_{L/X,x}
$
and the symmetric isomorphism is given by $q_H$, where one must
put in the twist to account for the scaling of $q_H$.

Differentiating again,
 we obtain
$$
<\lx 0' , \lh 0''>
= -<\lx 0''',\lh 0> -2<\lx 0 '' ,\lh 0 '>.
$$
The first term may be interpreted as
$w\intprod II_{X^*,H}(W,W)$ where
$W\in T_HX^*$ corresponds to $\lh 0'$ and $w\in N^*_HX^*$
corresponds to $\lx 0'$.
The second term can be interpreted as
$-n_H\intprod\pii_{X,x}(w,w,w)$
where $\pii_{X,x}$ is the cubic form (see \cite{Lan}).  
The third term can be interpreted as
$-2W\intprod II_{X,x}(w,w)$ where here we consider
$W\in N^*_xX$ and $w\in T_xX$ via the identifications
described above.  Polarizing, we have our complete inversion
formula:

\proclaim{Theorem 3.9}Let $H\in X^*_{sm}$ be  such
that $L_H\cap X_{sing}=\emptyset$. Let $x\in L_H$ be any point.
Then $II_{X^*,H}$ is determined by 
$II_{X,x}$ and $\pii_{X,x}$. More precisely,
given $y\in N^*_HX^*, v,w\in T_HX^*$, one has
the following formula defined up to twists and scales:
$$
\align 
y \intprod II_{X^*,H}(v,w)
 =& n_H\intprod \pii _{X,x}(  y,v_1,w_1)
+ n_H\intprod II _{X,x}(v_1,w_1)\tag 3.10\\
& 
 + v_2\intprod II _{X,x}(  y,w_1)
+w_2\intprod II _{X,x}(  y,v_1) \endalign
$$
All vector bundles in (3.10) are pulled back to $\ci$.
$n_H\in N^*_xX$ denotes a vector representing $H$.
$v_1,w_1\in T_HX^*/T_HL_x$ denotes the projection
of $v,w$ to
the quotient, and $v_2,w_2\in T_HL_x$
denotes their projection to the subspace.
All vectors appearing on the right hand
side of (3.10) are well defined elements up to twists and scales
of the appropriate spaces via repeated use of (3.5).
\endproclaim

The ambiguity  usually
occuring with the cubic form is 
eliminated by lifting to $\ci$ and evaluating
it on $n_H$ paired with at least one element of the singular
locus of $q_H$.

One can get a completely well defined formula using
the twists, but the resulting system of quadrics of
course will be the same.

\subheading{3.11 Some questions}
Do the examples of 2.12  for $a \ge 3$ arise from interesting 
dual varieties? Of course,
for odd $a$ we must first double to get a symmetric system. $a=3$ would
yield a 46 dimensional variety in $\pp{77}$ with defect 6 and $a=4$ would be
a 78 dimensional variety in $\pp{135}$ with defect 8.

In the boundary case of Theorem 2.16, i.e. $l = m - r$ (with $l \ge 1$),
the only examples we know have $E = T_{\pp{2}}(-1)^{\oplus 2}$, 
$\Lambda^2(T_{\pp{4}}(-1))$ or $\co_{\pp{l}}^{r/2} \oplus
\co_{\pp{l}}^{r/2}(1)$. Are there any others? Smooth varieties with
$\tdim(X) = \tdim(X^*)$ provide examples, but assuming Hartshorne's
conjecture, Ein's classification 
theorem \cite{E, Theorem 4.5} along with the remarks at the end of 2.13
show that no further
examples arise in this way. On the other hand, without assuming
Hartshorne's conjecture, Ein's result shows that if $l \ge r \ge 4$ then
$\co_{\pp{l}}^{r/2} \oplus \co_{\pp{l}}^{r/2}(1)$ does not arise from
dual varieties.

Suppose $\co_{\pp{l}}^m \to E \to 0$ and $E \cong E^*(1)$. If 
$r \le 2l-4$ is $E$ free? (One can also ask the same question with the
more restrictive hypotheses that the isomorphism is symmetric or
skew symmetric). The evidence for this is that this is the line
below $(l,n) = (2,2)$ and $(4,6)$ corresponding to the non-free examples
$T_{\pp{2}}(-1)$ and $\Lambda^2(T_{\pp{4}}(-1))$.

\enddocument